%% file: main_arxiv.tex
\documentclass[a4paper, 
10pt
]{article}

\input{sections/preamble.tex}

\begin{document}

\title{An efficient computational framework for naval shape design and
optimization problems by means of data-driven reduced order modeling
techniques} 
\author[1]{Nicola~Demo\footnote{nicola.demo@sissa.it}}
\author[1,2]{Giulio~Ortali\footnote{giulio.ortali@sissa.it}}
\author[3]{Gianluca~Gustin\footnote{gianluca.gustin@fincantieri.it}}
\author[1]{Gianluigi~Rozza\footnote{gianluigi.rozza@sissa.it}}
\author[3]{Gianpiero~Lavini\footnote{gianpiero.lavini@fincantieri.it}}
\affil[1]{Mathematics Area, mathLab, SISSA, International School of Advanced Studies, Trieste, Italy}
\affil[2]{Department of Mathematical Sciences, Politecnico di Torino, Italy}
\affil[3]{Fincantieri - Divisione Navi Mercantili e Passeggeri, FINCANTIERI SpA, Trieste, Italy}

\maketitle

\begin{abstract}
This contribution describes the implementation of a data--driven shape
	optimization pipeline in a naval architecture application. We adopt
	reduced order models (ROMs) in order to improve the efficiency of the
	overall optimization, keeping a modular and equation-free nature to
	target the industrial demand.  We applied the above mentioned pipeline
	to a realistic cruise ship in order to reduce the total drag. We begin
	by defining the design space, generated by deforming an initial shape
	in a parametric way using free form deformation (FFD). The evaluation
	of the performance of each new hull is determined by simulating the
	flux via finite volume discretization of a two-phase (water and air)
	fluid. Since the fluid dynamics model can result very expensive ---
	especially dealing with complex industrial geometries --- we propose also
	a dynamic mode decomposition (DMD) enhancement to reduce the
	computational cost of a single numerical simulation. The real--time
	computation is finally achieved by means of proper orthogonal
	decomposition with Gaussian process regression (POD-GPR) technique.
	Thanks to the quick approximation, a genetic optimization algorithm
	becomes feasible to converge towards the optimal shape.
\end{abstract}

\input{sections/intro.tex}

\input{sections/pipeline.tex}

\input{sections/results_arxiv.tex}

\input{sections/conclusion.tex}

\section*{Acknowledgements}
We thank Prof. Claudio Canuto for his constant support.

This work was partially supported by European Union Funding for Research and
	Innovation --- Horizon 2020 Program --- in the framework of European
	Research Council Executive Agency: H2020 ERC CoG 2015 AROMA--CFD
	project 681447 ''Advanced Reduced Order Methods with Applications in
	Computational Fluid Dynamics`` P.I. Gianluigi Rozza. The work was also
	supported by INdAM-GNCS: Istituto Nazionale di Alta Matematica --
	Gruppo Nazionale di Calcolo Scientifico.

%
%

\bibliographystyle{abbrv}      

\bibliography{gpr-biblio}   

\end{document}

%% file: sections/preamble.tex
\usepackage{graphicx}
\usepackage{latexsym}
\usepackage{titlesec}
\usepackage{hyperref}
\usepackage{subfig}
\usepackage{pgfplots}
\usepackage{pgfplotstable}
\usepackage{mathtools}
\usepackage{multicol}
\usepackage{comment}
\usepackage{booktabs}
\pgfplotsset{compat=1.5}
\usepackage{amssymb}
\usepackage{url}
\usepackage{bm}
\usepackage{pdfsync}
\usepackage{float}
\usepackage{tabularx}
\usepackage{enumerate}
\usepackage{array}
\usepackage{xspace}
\usepackage{wasysym}
\usepackage{empheq}
\usepackage{siunitx}
\usepackage{xfrac}
\usepackage{authblk}

\usetikzlibrary{spy}

\newcommand{\mupar}{\ensuremath{\boldsymbol{\mu}}}
\newcommand{\domain}{\Omega}
\newcommand{\pdomain}{\Omega(\bm \mu)}
\newcommand{\R}{\mathbb{R}}

\titleformat*{\subsection}{\bfseries}
\titleformat*{\subsubsection}{\bfseries \itshape}

%% file: sections/intro.tex
\section{Introduction and motivations}
\label{sec:intro}

A shape optimization problem consists of finding the geometric configuration of
an object that maximizes the performance of such object. Due to the number and
the complexity of methods to integrate together --- i.e. a shape
parametrization algorithm, a numerical solver, an optimization procedure ---,
this task remains challenging even nowadays. One of the most common problems
is the computational cost required to solve the mathematical model, necessary
to predict the performance of the deformed object. Addressing complex
phenomena, even exploiting high-performance facilities, the total computational
load may make the procedure unfeasible, since the performance evaluation has
to be repeated for each new deformed configuration.

In this work, we extend the computational pipeline already presented
in~\cite{DemoTezzeleGustinLaviniRozza2018NAV}, using two different reduced
order modeling (ROM) approaches to address the high computational demand of
optimization problems based on partial differential equations (PDEs) in
parametric domains. The goal is obtaining the optimal shape of the input object
--- in our case, the naval hull of a cruise ship --- with a resonable demand of
computational resources. For different version of this shape optimization
pipeline, we suggest~\cite{DemoTezzeleMolaRozza2019MARINE} for a POD reduction to geometrical parametrization
and~\cite{tezzele2019marine,mola2019marine} for an additional parameter space
analysis by means of active subspace property. ROM provides a model
simplification, bartering a slightly increased error in the model output with a
remarkable reduction of the computational cost. The real--time response of such
models helps to accelerate the entire optimization process. Other similar
framework regarding ROM have been presented
in~\cite{rozza2018advances,tezzele2019mortech}.

In details, the two adopted ROM techniques are: {\it i}) the data--driven
proper orthogonal decomposition (POD) coupled with Gaussian process regression
(GPR) for the approximation of the solution manifold for the
parametric model, and {\it ii}) the dynamic mode decomposition (DMD) algorithm
to estimate the regime state of the transient fluid dynamics problem.
Exploiting these techniques, not only we need a limited number of high--fidelity
(and expensive) simulations but we are able even to reduce the computational
cost of the latter. The main advantage is that the optimization
procedure, which has to iterate towards the optimum, uses the reduced order
model to estimate the performance of any new deformed object in a very quickly manner. 
An additional value of the proposed framework is the complete modularity for
the data-driven nature of the ROM methods. In fact, they are based only on the
output of the system, without the necessity to know the governing equations or,
from a technical viewpoint, to access to the discrete operators of the problem.
We propose in this work an application on the shape optimization of a cruise
ship, but the pipeline can be easily modified to plug different algorithms or
software.
All these features make the framework especially suited for industry,
thanks to the huge speedup in optimization --- but also design --- contexts and
the natural capacity to be even coupled with commercial software.

The work is structured as following: in \autoref{sec:pipeline} we described in
the details how the components are combined together, going into the deeper
mathematical formulation of all of them in the next subsections. In particular:
\autoref{sec:ffd} will focus on the free-form deformation (FFD), the algorithm
used for the shape parametrization; \autoref{sec:fom} will introduce the
full-order model we adopt and its numerical solution using the finite volume
(FV) approach; \autoref{sec:dmd} and \autoref{sec:pod} will introduce the
algorithmic formulation of the DMD and POD-GPR techniques, respectively;
\autoref{sec:opt} will summarize the genetic algorithm (GA), i.e. the optimization
method we used. Finally, in \autoref{sec:results} we present the numerical
setting of the resistance minimization problem for a parametric cruise ship and
the results obtained by applying the described framework on it, before
proposing a conclusive comment and some future perspectives in
\autoref{sec:conclusion}.

%% file: sections/pipeline.tex
\newcommand{\M}{\mathcal{M}}
\newcommand{\hdim}{\mathcal{N}}
\newcommand{\highspace}{\R^{\hdim}}
\newcommand{\ldim}{N_\text{POD}}
\newcommand{\lowspace}{\R^{\ldim}}
\newcommand{\out}{\mathbf{y}}
\newcommand{\cc}{\mathbf{c}}

\section{The complete computational pipeline}
\label{sec:pipeline}
\begin{figure*}
\includegraphics[width=\textwidth]{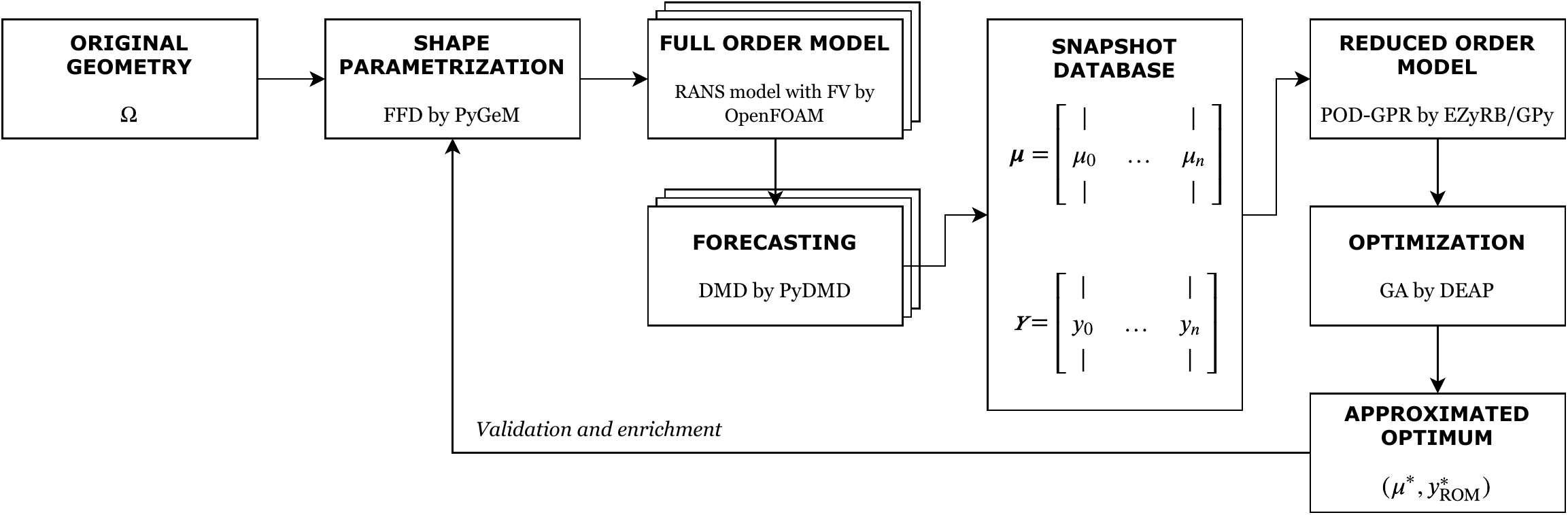}
\caption{The flowchart of the complete computational pipeline.\label{fig:scheme}}
\end{figure*}
This section focuses on the integration of all the components into a single
pipeline capable of optimizing an input object with a generic shape $\domain
\in \R^3$. We will provide details about the methodologies stack, specifying
the interfaces between methods in order to let the reader capable to understand
the workflow. The proposed framework can be however, thanks to the data-driven
feature, easily extended, replacing one or more tecnhniques, increasing
integrability of such pipeline.

As first ingredient, we need a map
$\M:\R^3 \to \R^3$ that,
depending on some numeric parameters, deforms the original domain such that
$\pdomain = \M(\domain, \mupar)$. Dealing with complex geometries, we chose the
free-form deformation (FFD)~\cite{sederberg} to deform the original object, because of its
capability to preserve continuity on the surface derivatives and to perform
global deformation even with few parameters. The parameters $\mupar \in P \subset
\R^P$, for this method, control the displacement of some points (along some
directions) belonging to a lattice of points around the object. This motion
produce a deformation in all the space embedded by the lattice.
Chosen the parameter space $P$, we sample this latter $N$ times to obtain the
set $\{\mupar_i\}_{i=1}^N$, and, using the FFD, the corresponding set of
deformed shapes $\{\domain(\mupar_i)\}_{i=1}^N$.\smallskip

The performance of all the samples have to be evaluated, using an accurate
numerical solver. In this case, since the analyzed problem is related to an
incompressible turbulent multiphase flow, we use the Reynolds-averaged
Navier--Stokes (RANS) equations with the volume of fluid (VOF) approach to
describe the mathematical model, and a finite volume (FV) discretization to
numerically solve it.
Such model requires, both for the complex geometry and the complexity of
equations, a not-negligible amount of computational resources. Even
if, as in this case, the number of these high-fidelity simulations is limited
to $N$, the overall load may result too big. We can gain additional speedup
exploiting the dynamic mode decomposition (DMD)~\cite{kutz} to predict the
regime state of the simulation. In our test, the time-dependent problem shows a
quasi-periodic behaviour, continuing to oscillate around the asymptotic
configurations. DMD catches this kind of patterns in the temporal evolution of
a system, allowing to easily make predictions with a good accuracy.
We can combine the two techniques, by computing the initial temporal snapshots
--- aka the output of interest of such system at a certain time --- with
the high-fidelity model, then feeding the DMD algorithms with the latter in order to
predict the regime snapshots. We define the snapshots $\out_i^k$ as the output
of interest of the parametric domain $\domain(\mupar_i)$ at time $k$: the
regime state $\out_i^{m+c}$ is then predicted collecting the snapshots
$\{\out_i^{j}\}_{j=0}^m$, for $i = 1, \dotsc, N$. 
It is important to specify that the computational grids built around the
objects $\domain(\mupar_i)$ are not enforced to share the same topology, or the
same number of degrees of freedom, but for the DMD is necessary that the grids
do not change during the temporal evolution of the system.
In this work we do not use the pressure and velocity fields as output of
interest, but directly to the distribution of total resistance (over the
surface of hull). Since the data-driven approach, this does not imply any
additional complexity. Our database contains thus the discrete distribution of
the total resistance for all the samples.\smallskip

After this step, we obtain a set of $N$ pairs composed by the input parameters
and the regime states, that is $\{(\mupar_i, \out_i^{m+c})\}_{i=1}^N$. In case
of output with different dimensions, we need to project the solution from the
FV discretized space to the original deformed geometry $\pdomain$. Being
originated by the FFD, all the geometries share the same topology. Assuming the
geometry $\domain$ is discretized in $\hdim$ degrees of freedom, the resulting
new pairs are defined as $(\mupar_i, \hat{\out}_i)$, with $\mupar_i \in P$ and
$\hat{\out}_i \in \R^\hdim$, for $i = 1, \dotsc, N$.\smallskip

Proper orthogonal decomposition (POD)~\cite{rozza} is now
involved to reduce the dimensionality of the snapshots. The outputs
$\hat{\out}_i \in \highspace$ are
projected onto the POD space, which typically has a very lower dimensions,
obtaining the reduced space representation $\cc_i \in \lowspace$ of the
original states. The input-output pairs are now $(\mupar_i, \cc_i)$ for $i = 1,
\dotsc, N$: assuming that a mapping $F:P \to \lowspace$ exists between input
and output such that $\cc = F(\mupar)$, we can exploit the collected outputs
to approximate the output itself for different parameter value using any
interpolation/regression method. In this contribution, we adopt a Gaussian
process regression (GPR)~\cite{quinonero2005unifying} to approximate the
input-output relation with a Bayesian approach. Other examples for the POD-GPR
coupling can be found in~\cite{guo2018reduced,OrtaliDemoRozzaCanuto2020}.
Finally, the low-dimensional output is projected back to the full-order space
to obtain the approximated solution. Combining the techniques, we are able to
build a reduced order model based only on the system output capable to provide
an approximation of the output $\out^j_\text{ROM}$ for untried parameters
$\mupar_j$ in real-time. In our test, we remember we use the resistance
distribution as output of interest.\smallskip

The optimization procedure is then applied over the reduced order model, by
computing the objective function on the state predicted using POD-GPR. Thanks
to the negligible time required for the performance evaluation of a new shape,
we can explore the parameter space with a genetic algorithm
(GA)~\cite{holland1973genetic} to converge to the optimal shape. The quanity to
minimize, in our numerical experiments, is the total resistance, that is
nothing but the intergral of corresponding field.  The objective function
relies hence on the previously mentioned methods, since to compute it we need
to project the POD-GPR approximation over the new shape obteined by FFD.
At the end, we get the optimal parameter $\mupar^*$ and
correspondent (approximated) output $\out^*_\text{ROM}$.  Such parameter can be
used to restart the pipeline, performing the morphing over the geometry then
testing it by using the high-fidelity solver for the validation of the result. Non
only: this latter simulation can be further exploited by adding it to the
snapshots database, resulting in an iterative process where the approximated
output is used for the reduced order model, enriching in this way the accuracy
of the model itself.

\input{sections/ffd.tex}

\input{sections/fom.tex}

\input{sections/dmd.tex}

\input{sections/pod.tex}

\input{sections/opt.tex}

%% file: sections/ffd.tex
\subsection{Free-form deformation for shape parametrization}
\label{sec:ffd}

Free-form deformation (FFD) is a geometric tool, extensively employed in
computer graphics, used to deform a rigid object based on the movement of some
predefined control points. Introduced in
\cite{sederberg}, it has seen various improvements over
the years. The reader can refer for example for a
more recent review~\cite{intr1} and~\cite{LassilaRozza2010a}
and~\cite{SalmoiraghiScardigliTelibRozza2018} for a coupling with ROM
techniques.
The main idea behind FFD is to define a regular lattice of points around the
object (or part of it) and manipulate the whole embedded space by moving some
of those control points.  Mathematically, this is obtained by mapping the
physical space enclosed by the lattice to a unit cube $D = [0, 1]^d$ by using an
invertible map $\psi:\R^d \to D$.

Inside the unit cube we define a cubic lattice of control points, with $L$,$M$
and $N$ points respectively in $x$,$y$ and $z$ directions:
\begin{equation}
P^0_{l,m,n}=\begin{pmatrix}
l/L \\
m/M \\
n/N
\end{pmatrix} \in D,
\end{equation}
where $l=0,\dotsc,L$, $m=0,\dotsc,M$ and $n=0,\dotsc,N$ .
We move these points by adding a motion $\mu_{l,m,n}$
such that:
\begin{equation}
P_{l,m,n} = P^0_{l,m,n} + \mupar_{l,m,n}.
\end{equation}
The parametric map $\hat{T}:D \to D$ that performs the deformation of reference space is then defined by:
\begin{equation}
\hat{T}(s,t,p; \mupar) = \sum_{l=0}^L \sum_{m=0}^M \sum_{n=0}^N 
b_{l}^{L}(s)
b_{m}^{M}(t)
b_{n}^{N}(p)
P_{l,m,n},
\end{equation}
where:
\begin{equation}
\begin{array}{l}
b_{l}^{L} (s) = \binom{L}{l} (1-s)^{(L-l)}s^l,\\
b_{m}^{M} (t) = \binom{M}{m} (1-t)^{(M-m)}t^m,\\
b_{n}^{N} (p) = \binom{N}{n} (1-p)^{(N-n)}p^n.\\
\end{array}
\end{equation}
The FFD map $T:\R^3 \to \R^3$ is then composed as it follows:
\begin{equation}
T(\cdot; \mupar) = (\psi^{-1} \circ \hat{T} \circ \psi)(\cdot; \mupar).
\end{equation}
We applied the FFD algorithm directly to input object using the open source
Python package called PyGeM~\cite{pygem}.

%% file: sections/fom.tex
\subsection{Finite volume for high-fidelity database}
\label{sec:fom}
We now discuss the full order model (FOM), which generates what we call the
high fidelity solution.
The Reynolds--averaged Navier--Stokes (RANS) equations \\model the turbulent
incompressible flow around the naval hull, while for
the modeling of the two different phases --- e.g. water and air --- we adopt the
volume of fluid (VOF) technique~\cite{vof}.
The equations governing our system are then:
\begin{equation}
\begin{cases}
\frac{\partial \bar{u}}{\partial t} + (\bar{u} \cdot \nabla)\bar{u}
+\frac{1}{\rho} \nabla \bar{p} -\nabla \cdot \nu \nabla \bar{u} - \nabla \cdot
(\tilde{u} \otimes \tilde{u})=0,\\
\nabla \cdot \bar{u}=0,\\
\frac{\partial \alpha}{\partial t} + \nabla \cdot (\bar{u} \alpha)=0,\\
\end{cases}
\label{ns}
\end{equation}
where $\bar{u}$ and $\tilde{u}$ refer the mean and fluctuating velocity after
the RANS decomposition, $\bar{p}$ denote the (mean) pressure, $\rho$ is the
density, $\nu$ the kinematic viscosity and $\alpha$ is the discontinuous
variable belonging to interval $[0, 1]$ representing the fraction of the second
flow in the infinitesimal volume.

The first two equations are the continuity and momentum conservation, while the third one represent the transport equation for the VOF variable $\alpha$.
The Reynolds stresses tensor $\tilde{u} \otimes \tilde{u}$ can be modeled by
adding additional equations in order to close the system: in this work, we use
the $\text{SST} k-\omega$ turbulence model~\cite{menter}.  For the
multiphase nature of the flow, the density $\rho$ and the kinematic viscosity
$\nu$ are defined using an algebraic formula expressing them as a convex
combination of the corresponding properties of the two flows:
\begin{equation}
\begin{array}{l}
\rho = \alpha \rho_1 + (1-\alpha) \rho_2,\\
\nu = \alpha \nu_1 + (1-\alpha) \nu_2.
\end{array}
\label{density}
\end{equation}

To solve such problem, we apply the finite volume (FV) approach. We adopted a
$1^{st}$ order implicit Euler scheme for the temporal discretization, while for
the spatial scheme we apply the linear upwind one. Regarding the software, the
simulation is carried out using the C++ library OpenFOAM~\cite{openfoam}.

%% file: sections/dmd.tex
\subsubsection{Dynamic mode decomposition for regime state prediction}
\label{sec:dmd}

Dynamic mode decomposition (DMD) is a data-driven ROM technique that
approximates the evolution of a complex dynamical system as the combination of
few features linearly evolving in time~\cite{schmid2010dynamic,kutz}. The basic idea is to provide a
low-dimensional approximation of the Koopman operator~\cite{koopman} based on
few temporarily equispaced snapshots of the studied system.
DMD assumes the evolution of the latter can be expressed as:
\begin{equation}
\out_{k+1} = \mathbf{A} \out_k,
\label{yax}
\end{equation}
where $\out_{k+1} \in \highspace$ and $\out_k \in \highspace$ are two snapshots
at the time $t = k$ and $t = k+1$, respectively, while $\mathbf{A}$ refers to
a discrete linear operator. A least-square approach can be used to calculate
this operator. After collecting a set of snapshots defined as
$\{\out_{t_0+k\Delta t}\}_{k=0}^M$, we can arrange them into two matrices
$\mathbf{Y} = \begin{bmatrix} \out_0 & \dotsc & \out_{M-1} \end{bmatrix},
	\dot{\mathbf{Y}} = \begin{bmatrix} \out_1 & \dotsc & \out_{M}
	\end{bmatrix} \in \R^{\hdim\times M}$ such that the correspondent
	columns of the two matrices represent two sequential snapshots.

We can now minimizing the error $\|\mathbf{A}\mathbf{Y} - \dot{\mathbf{Y}}\|_F$
by the following matrix multiplication $\mathbf{A} = \dot{\mathbf{Y}}
\mathbf{Y}^\dagger$, where the symbol $^\dagger$ indicates the Moore-Penrose pseudoinverse.
While we can already use the operator $\mathbf{A}$ to analyze the system, in
practice because of its considerable dimension and the difficulties that would
arise in order to obtain it numerically. DMD uses then the singular value
decomposition (SVD) to compute the reduced space onto which projecting the
operator. Formally
\begin{equation}
\mathbf{Y=U \Sigma V}^T,
\end{equation}
where $\mathbf{U} \in \mathbb{R}^{\hdim\times M}$, $\mathbf{\Sigma} \in
\mathbb{R}^{M\times M}$ and $\mathbf{V} \in \mathbb{R}^{M\times M}$. The left
singular vectors (the columns of $\mathbf{U}$) span the optimal low-dimensional
space, allowing us to project the operator $\mathbf{A}$ onto it:
\begin{equation}
	\mathbf{\tilde{A} = U}^T \mathbf{A U = U}^T \mathbf{\dot{Y} V \Sigma}^{-1}
\end{equation}
to compute the reduced operator. The interesting feature is that the
eigenvalues of $\mathbf{\tilde{A}}$ are equal to the non-zero ones of the high
dimensional operator $\mathbf{A}$, and also the eigenvectors of the two operators are related each other~\cite{tu}. In particular:
\begin{equation}
	\mathbf{\Phi = \dot{Y} V \Sigma}^{-1} \mathbf{W},
\end{equation}
where $\Phi$ is the matrix containing the $\mathbf{A}$ eigenvectors, the
so-called DMD modes, and $\mathbf{W}$ is the matrix of $\tilde{\mathbf{A}}$
eigenvectors.
Defining $\mathbf{\Lambda}$ as the diagonal matrix of eigenvalues, we have:
\begin{equation}
\mathbf{Y} \approx \mathbf{A} \mathbf{X} = \mathbf{\Phi} \mathbf{\Lambda} \mathbf{\Phi}^\dagger \mathbf{X}
\end{equation}
that implies that any snapshots can be approximated computing $\out_k =
\mathbf{\Phi} \mathbf{\Lambda}^k \mathbf{\Phi}^\dagger \out_0$.

We apply the DMD on the snapshots coming from the full-order model (discussed
in \autoref{sec:fom}) in order to perform fewer temporal iterations using the
high-fidelity solver, and predict the output we are interested to analyze in
order to gain an additional considerable speedup. The results are obtained
using PyDMD~\cite{pydmd}, a Python package that implements the most common
version of DMD.

%% file: sections/pod.tex
\subsection{Reduced order model exploiting proper orthogonal decomposition}
\label{sec:pod}

Reduced basis (RB) is a ROM method that
approximates the solution manifold of a parametric
problem using a low number of basis functions that
form what we call the reduced
basis~\cite{morhandbook,rozza}. In this community,
proper orthogonal decomposition (POD) is a widespread
technique~\cite{struct,StaRo2018} since its capability to provide
orthogonal basis that have an energetically hierarchy. While a possible
approach for turbulent flows involving projection-based ROM is available
in~\cite{HiStaMoRo2019}, we prefer the data-driven approach for the higher
integrability in many industrial workflows. POD needs as
input a matrix containing samples of the solution
manifold. We
define $\hdim$ the number of degrees of freedom of our numerical model and
$\out \in \highspace$ its solution for a generic parameter $\mupar$. Thus, the
snapshots matrix $\mathbf{Y} \in \R^{\hdim\times n}$ is defined as:
\begin{equation}
\mathbf{Y} = 
\left[
  \begin{array}{cccc}
    \vrule & \vrule & & \vrule\\
    \out_{1} &\out_{2} & \ldots & \out_{n} \\
    \vrule & \vrule & & \vrule 
  \end{array}
\right].
\end{equation}
The POD basis is defined as basis that maximizes the similarity (as measured by
the square of the scalar product) between the snapshots matrix and its
elements, under the constraint of orthonormality. Formally, the POD basis $\{\psi_i\}_{i=0}^l$ of dimension $l$ is defined as:
\begin{equation}
\max_{\psi_1,\dotsc\psi_l} \sum_{i=1}^l \sum_{j=1}^n
| \langle\mathbf{y}_j,\psi_i\rangle_{\highspace} |^2
\label{pod_basis}
\end{equation}
such that
$\langle\tilde{\psi_i},\tilde{\psi_i}\rangle_{\highspace}=\delta_{i,j}$, for $1
\leq i,j \leq l$.
Singular value decomposition (SVD) is a method that computes the POD basis~\cite{volkwein-rom} 
by decomposing the snapshots matrix: 
\begin{equation}
\mathbf{Y} = \mathbf{U} \mathbf{\Sigma} \mathbf{V}^*,
\end{equation}
where matrices $\mathbf{U} \in \R^{\hdim\times n}$ and $\mathbf{V} \in \R^{n
\times n}$ are unitary while $\mathbf{\Sigma} \in \R^{n\times n}$ is diagonal.
In particular, the columns of $\mathbf{U}$ are POD basis.
We project the original snapshots onto the POD space to have a low-dimensional representation. In matrix form:
\begin{equation}
\mathbf{C} = \mathbf{U}^T \mathbf{Y}, \quad \mathbf{c} \in \R^{n\times n}.
\end{equation}
The columns of $\mathbf{C}$ are the modal coefficients $\mathbf{c}_i \in \R^n$.

\newcommand{\Xtrain}{\mupar}
\newcommand{\Xtest}{\bar{\mupar}}
\newcommand{\ytrain}{\mathbf{c}}
\newcommand{\ytest}{\bar{\mathbf{c}}}
We can now exploit this reduced space in order to build a probabilistic
response surface using the Gaussian process regression
(GPR)~\cite{quinonero2005unifying}. In particular,
assuming there is a natural relation $F:P \to \R^n$ between our geometric
parameters $\mupar$ and the low-dimensional output $\mathbf{c}$ such that
$\mathbf{c} = F(\mupar)$, we try to approximate it with a multivariate Gaussian
distribution. We define:
\begin{equation}
f(\mupar) \sim \text{GP}(\mathcal{M}(\mupar), \mathcal{K}(\mupar, \mupar)),
\label{eq:gp}
\end{equation}
where $\mathcal{M}$ refers to the mean of the distribution and $\mathcal{K}$ to
its covariance. There are many possible choices for the covariance function
$\mathcal{K}:P \times P \to \R$, in our case we use the the squared exponential
one defined as $\mathcal{K}_{SE}(\mupar_i, \mupar_j) = \sigma^2 \exp(-\frac{1}{2}\|\mupar_i -
\mupar_j\|^2)$.
The prior joint Gaussian distribution for the outputs $\mathbf{c}$ results then
\begin{equation}
\mathbf{c}|\mathbf{\mupar}\sim~\mathcal{N}(0,\mathcal{K}(\Xtrain, \Xtrain)).
\end{equation}
For sake of simplicity we assume that the GP has mean equal to zero: the entire
process results defined only by the covariance function. In order to specify
the GP for our dataset, we need to maximize the marginal likelihood varying the
hyper-parameters of the covariance function, in this case only the $\sigma$.
Once obtained the output distribution, we can just sample it at the test
parameters to predict the output --- which, we remember, is the low-dimensional
snapshot --- by exploiting the joint distribution:
\begin{equation}
\ytest | \Xtest, \Xtrain, \ytrain \sim \mathcal{N} (\mathbf{m}, \mathbf{C})
\end{equation}
with 
\begin{equation}
\begin{array}{ll}
\mathbf{m} &= \mathcal{K}(\Xtest,
\Xtrain)\mathcal{K}(\Xtrain,\Xtrain)^{-1}\ytrain,\\
\mathbf{C} &=
\mathcal{K}(\Xtest, \Xtest) -
\mathcal{K}(\Xtest,\Xtrain)\mathcal{K}(\Xtrain,\Xtrain)^{-1}\mathcal{K}(\Xtrain,\Xtest),
\end{array}
\end{equation}
where $\mupar$ and $\bar{\mupar}$ refer to the input parameters and the test
parameters, and where $\mathbf{c}$ and $\bar{\mathbf{c}}$ are the corresponding train
and test output.

We compute the modal coefficients of all (untested) new parameters. To
approximate the high-dimensional snapshots we need just to back map the modal
coefficients to the original space. In matricial form:
\begin{equation}
\bar\out = \mathbf{U}\bar{\mathbf{c}}
\end{equation}

An additional gain of such method is the complete division between two
computational phases often called {\it offline} and {\it online} steps We can easily
note that, to collect the input snapshots, we initially need to compute several
snapshots using the chosen high-fidelity model. This is the most expensive
part, and usually is carried out on powerful machines. This offline step is
fortunately independent from the online one, where actually the snapshots are
combined to span the reduced space and approximate the new reduced snapshots.
Since this latter can be easily performed on standard laptops, the
computational splitting in two steps allows also to efficiently exploit all
different resources.

For the implementation, we developed and released this verion of data-driven
POD in the numerical open source package EZyRB~\cite{ezyrb}, exploiting the
library GPy~\cite{gpy2014} for the GPR step.

%% file: sections/opt.tex
\subsection{Genetic algorithm for global optimization}
\label{sec:opt}

Genetic algorithms (GA) denote in literature the family of computational methods
that are inspired by Darwin’s theory of evolution.  In an optimization context,
emulating the natural behaviour of living beings, this methodology gained
popularity due to its easy application and the capability to not get blocked in local
minima. The algorithm was initially proposed by Holland in~\cite{holland1973genetic,holland1992adaptation,back1996evolutionary} and it is based on few fundamental steps: selection, mutation and mate.
We consider any sample of the parameter domain as an individual $\mupar_i \in P
\subset \R^P$ with $P$ chromosomes. The fitness of the individuals is quantify
by a scalar objective function $f:P \to \R$.
We define the initial population $\mathbf{M}^0 = \{\mupar_i\}_{i=1}^{N_0}$
composed by $N_0$ individuals that are randomly created within the parameter
space. The corresponding fitnesses are compute and the evolutive process of
individuals starts.

The first step is the {\it selection} of the best individuals in the population.
Intuitively, the basic approach results chosing the $N$ individuals that have
the highest fitnesses, but for large population, or simply to reinforce the
stocastich component of the method, a propabilistic selection can be performed.
The selected individuals are often refered as the offspring that will breed the
future generation.

We are now ready to reproduce the random evolution of such individuals. This is
done in the {\it mutation} and {\it mate} steps. In the mutation, chromosomes
of the individuals can change, partially or entirelly, in order to create the
new individuals. Several approaches are available for the mutation, but usually 
they are based to a mutation probability to reproduce the aleatory nature of
evolution.  In the mate step, individuals are coupled into pairs and, still
randomly, the chromosomes of the parent individuals are combined to originate
the two children. In particular, the mate step emulates the
reproduction step, and for this reason can be usually called also cross-over. 

The population is now composed by the new (mated and mutated) individuals.
Iterating this process, the population will converge toward the optimal
individual, but depending on the shape of fitness function it may requires many
generations to converge. 

For the numerical experiments, we use the GA implementation provided by the
DEAP~\cite{deap} package, an open source library for evolutionary algorithms.

%% file: sections/results_arxiv.tex
\section{Numerical results: a cruise ship shape optimization}
\label{sec:results}

In this section, we will present the results obtained by applying the described
computational pipeline to optimize the shape of a cruise ship.
We maintain the same structure
of the previous section, discussing the intermediate results for any mentioned
technique.
\begin{figure}[h]
\centering
\includegraphics[width=.85\textwidth]{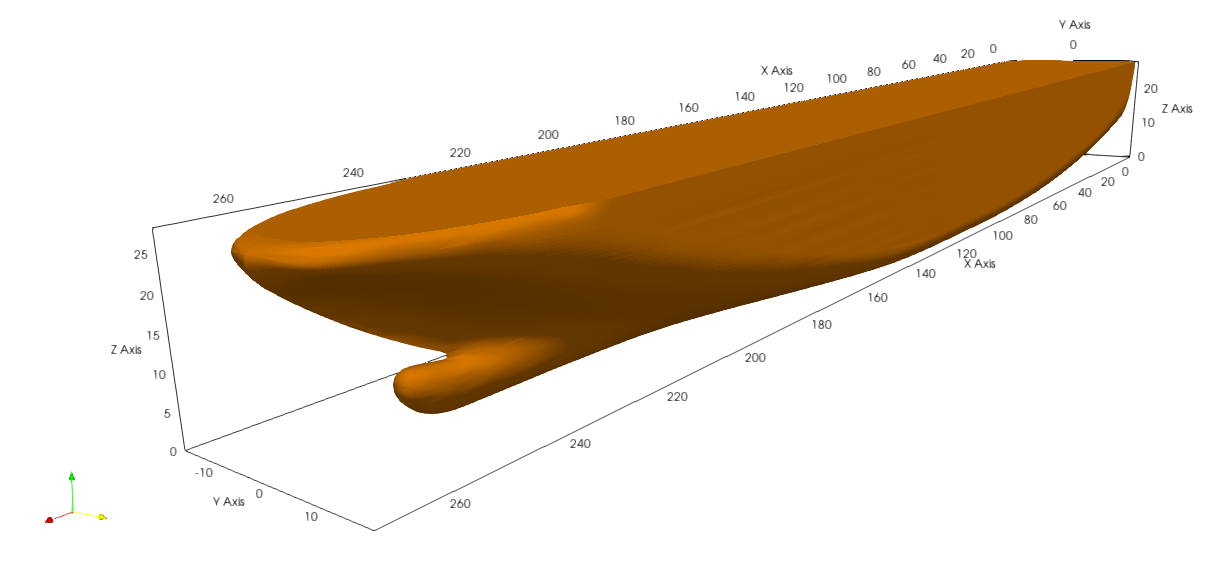}
\caption{The undeformed hull.}
\label{boat}
\end{figure}

\begin{figure}[h]
\centering
\includegraphics[width=.85\textwidth]{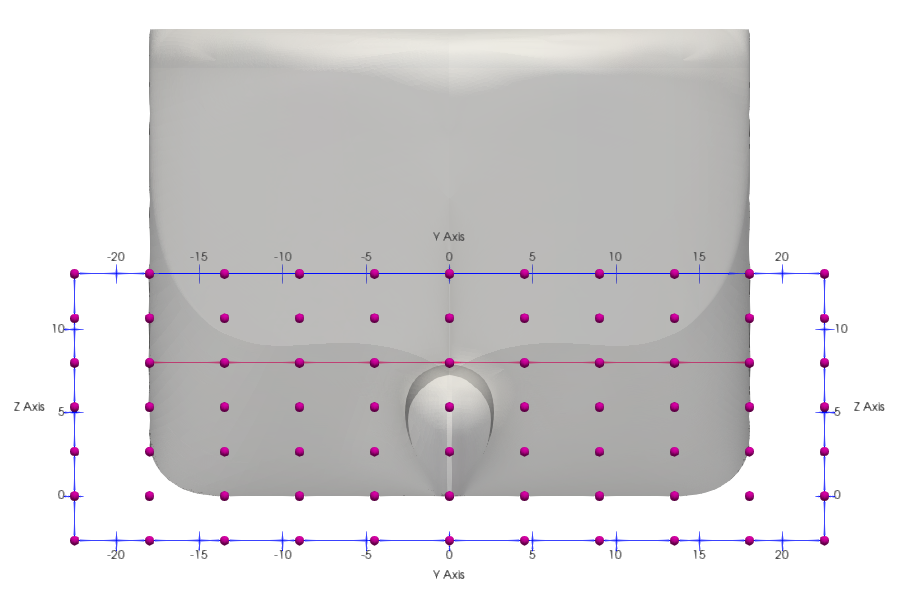}
\includegraphics[width=.85\textwidth]{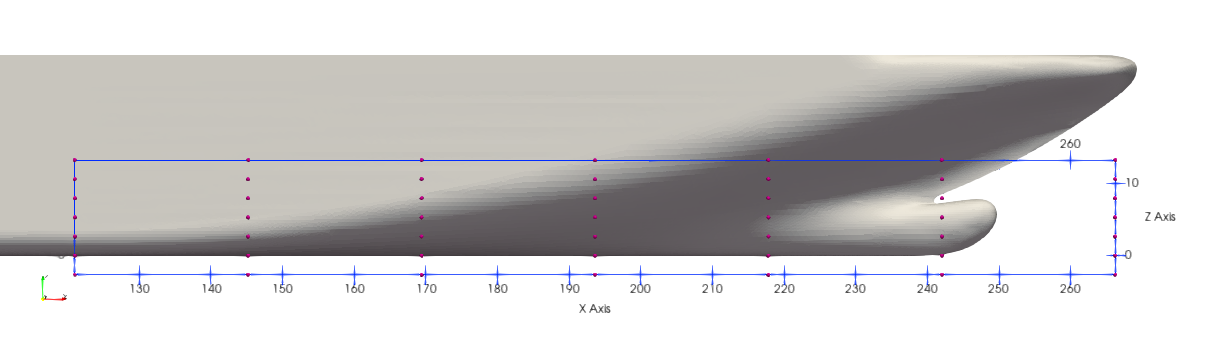}
\caption{$x$- and $y$-normal view of the set $D$ (in blue) and the lattice of points $P^0_{l,m,n}$ (in red) over the undeformed hull.}
\label{grid2}
\end{figure}

\paragraph{Free-form deformation}
We set the domain $D$, aka the space enclosed by the lattice of FFD control points, in order to deform only the immersed part of the hull, in the proximity of the bow.
The lattice is illustrated in \autoref{grid2}, and we can see that it is
positioned, in $x$ direction, on sections $10, 12, 14, 16, 18, 20$ and
22\footnote{In naval architecture a boat is divided, no matter the size, in $20$
chunks, generated by $21$ equally spaced cuts obtained with planes perpendicular
to the $x$-axis}. For $z$ direction, the points are displaced around the
waterline, while along $y$ axis the points are positioned for the entire width
of the ship. In total, $539$ FFD points are used.

\begin{figure}[h]
\centering
\includegraphics[width=.85\textwidth]{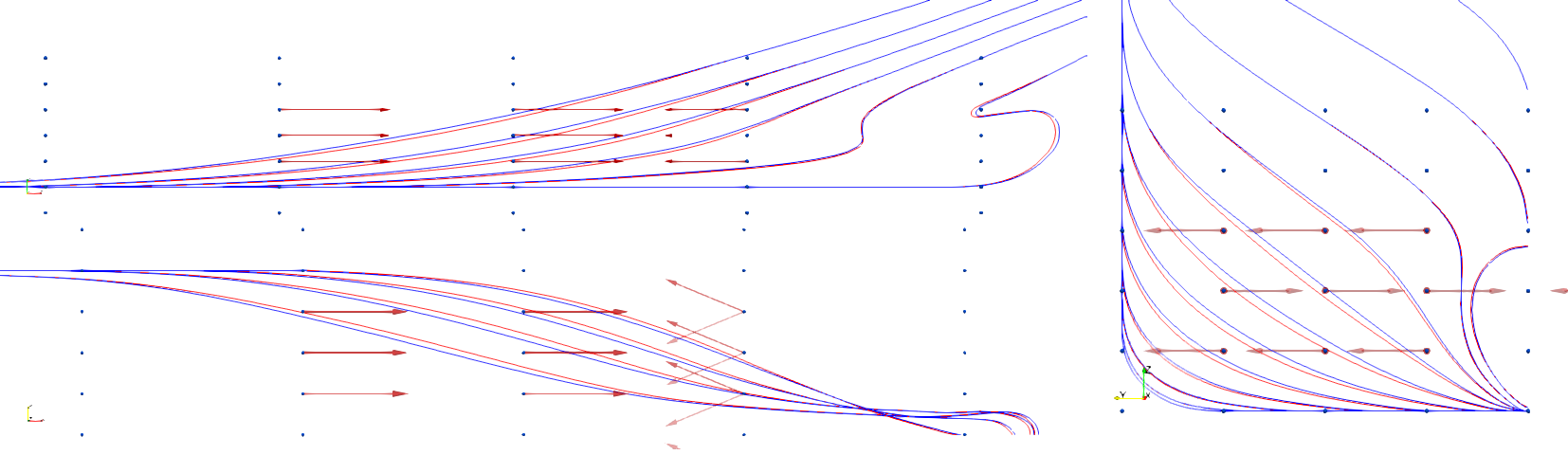}
\caption{Example of shape morphing with $\mupar~=~[0.08, 0.08, -0.06, 0.08, -0.08, 0.08]$}
\label{fig:ffd_ex}
\end{figure}

Concerning the motion of such points only part of the points in the lattice are
displaced: we use $6$ parameters to control the movements along $x$ (the first
three parameters) and along $y$. An example of this motion is sketched in
\autoref{fig:ffd_ex}, where red arrows refers to control points movements. The layers corresponding to sections $10, 12,
20$ and $22$ remain fixed, together with the two upper and lower layers, the
two far left and the two far-right layers and, finally, the layer over the
longitudinal symmetry plane. Except for this last one, that is kept fixed to
maintain symmetry, the other layers are kept fixed in order to achieve the
continuity and smoothness of the shapes, required especially in the $x$
direction where the deformation must link in a smooth way to the rest of the
boat.

The parameter range have been chosen in order to avoid a high decrease of the
hull volume and, at the same time, explores a large variety of new shapes. In
details, we have a tolerance of the $1\permil$ for the volume constraint. With
a trial and error approach we define the parameter ranges, obtaining a
parameter space  that is $P = [-0.08, 0.08]^6$ (the dimension of such space is
the number of parameters, $6$ in this test).
We underline that the parameters refer to the motion normalized for the $D$
length along the corresponding direction.

We create a set of $100$ samples taking with uniform distribution on the
parametric space. These are the input parameters of the high-fidelity database
required for ROM.

\begin{figure}[!t]
\centering
\includegraphics[width=.85\textwidth]{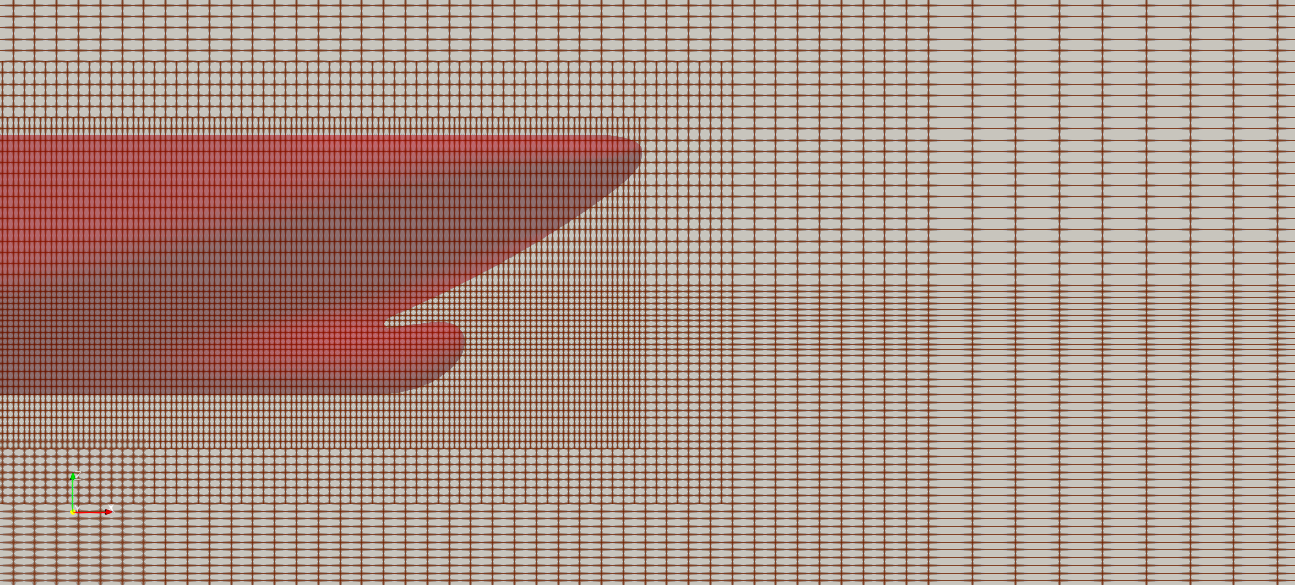}
\caption{The refined computational grid}
\label{fig:mesh}
\end{figure}

\paragraph{Finite volume discretization}
We simulate the flow pasting around the ship using the FV method, computing for
each deformed object the distribution of the total resistance over the hull.
The simulations are run on model scale ($1:25$).
The computational grid (defined in $[-39,24]\times[-29,0]\times[-24,6]$) is
built from scratch around all the deformed hulls, enforcing the mesh
quality. The computational grid counts $\approx\num{1.5e6}$ cells. To the VOF
model, we need an extra refinement around the waterline in order to avoid a
diffusive behavior of the fraction variable $\alpha$, which is discontinuous. A
region of the computation grid is reported in \autoref{fig:mesh} for
demonstrative purpose.
The numerical schemes adopted are mentioned in \autoref{sec:fom}, and we report
in \autoref{tab:fom} the main physical quantities we fix in our setting. The
Reynolds number is near to $\num{2e7}$.
\begin{table}[b]
\caption{Summary of the numerical variable for VOF-RANS model\label{tab:fom}}
\begin{tabular}{rl}
	inlet velocity $u$ & $\SI{2.26336}{\meter\per\second}$\\
	water density $\rho_W$ & $\SI{1.09e-6}{\meter^2\per\second}$\\
	air density $\rho_A$ & $\SI{1.09e-5}{\meter^2\per\second}$\\
	water kinematic viscosity $\nu_W$ & $\SI{998.8}{\kilogram\per\meter^3}$\\
	air kinematic viscosity $\nu_A$ & $\SI{1}{\kilogram\per\meter^3}$\\
	dissipation rate $\omega$ & $\SI{70.49721}{\second^{-1}}$\\
	turbulent kinetic energy $k$ & $\SI{7.6841e-4}{\meter^2\per\second^2}$\\
\end{tabular}
\end{table}
The integration in time is carried out for $t \in [0, 40]~\si{\second}$, with an initial
step of $\Delta t = \SI{1e-3}{\second}$ and an adjustable time-stepping governed by the
Courant number (we impose it to be lower than $5$). We clarify that, even if
the time stepping is not fixed, we save the equispaced temporal snapshots of
the system in order to feed the DMD algorithm.
In this work, we are interested to the total resistance of the ship: after
computing the pressure, velocity and fraction variable unknowns (from the
VOF-RANS model), we can exploit them in order to calculate the resistance distribution
(both the viscous and the friction terms) over the hull surface.
Regarding the computing time, on a parallel architecture with $40$ processes,
the simulation lasts approximately $8$ hours.

\paragraph{Dynamic mode decomposition}
We applied DMD on the results of the FOM. It is important to specify that we
fit a DMD model for each geometric deformation, as a sort of post-processing on
the output. We train the model using the snapshots $\{\out_i\}_{i=0}^{40}$ for
$i = t_\text{dmd}^0+i\Delta t_\text{dmd}$, with $t_\text{dmd}^0 = \SI{20}{s}, \Delta
t_\text{dmd} = \SI{0.5}{s}$. The first $\SI{20}{s}$ of the simulation are discarded since
they are not particularly meaningful for the boundary conditions propagation.
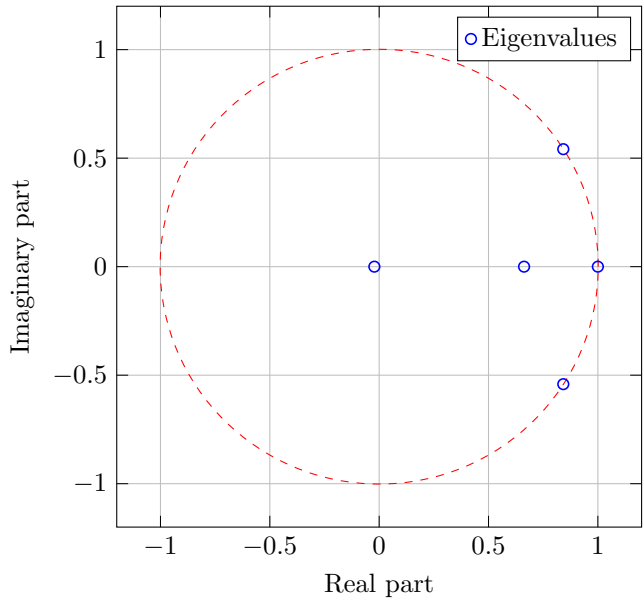
\begin{figure}[t]
\centering
\input{plots_arxiv/eigs.tikz}
\caption{First 5 DMD eigenvalues for a single simulation.}
\label{fig:eigs}
\end{figure}
In this contribution, we analyze the DMD operator from a spectral perspectives.
\autoref{fig:eigs} reports in fact the eigenvalues (computed for a single
simulation) after projecting the output onto a POD space of dimension $5$. The
position of eigenvalues in the complex plane provide information about the
dynamics of all the DMD modes. In particular, the imaginary part is related to
frequency, while the distance between them and the unit circle is related to
the growth-rate. We can neglect the dumped modes (the two eigenvalues inside
the circle) since their contribution is useless for future dynamics and focus
on the remaining ones: two modes present a stable oscillatory trend, that actually
catch the asymptotic oscillations of the FOM, and the last one ($1 + 0i$) is
practically constant. 
We isolate the contribution of only this latter mode, assuming it represents
the regime state to which the FOM converges, using it as final output.  In our
setting, having built the computational grid for all the deformed ships from
scratch, we need as last step to project the resistance distribution over the
initial geometry $\pdomain$, in order to ensure same dimensionality for all the
outputs. In our case we use a closest neighbors interpolation.
Thanks to the application of DMD, we can perform fewer time iterations in the
full-order model: in this case we can reduce the simulated seconds from $60$ to
$40$, approximating the regime state with DMD. This of course implies a
reduction of $\sfrac{1}{3}$ of the overall time required to run all the
simulations.

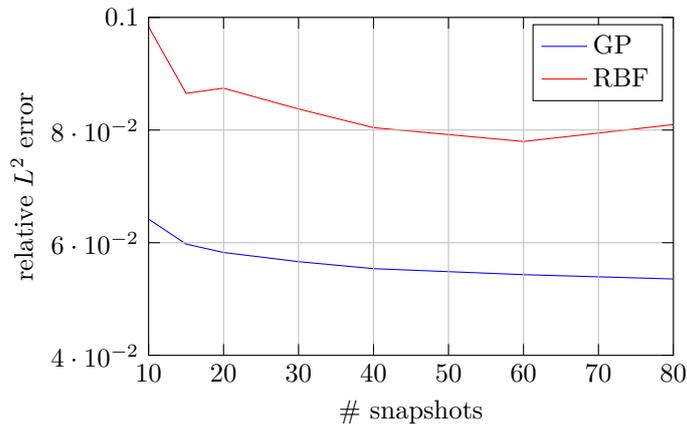
\begin{figure}[t]
\centering
\input{plots_arxiv/err2snapshots.tikz}
\caption{Sensitivity analysis on the accuracy of POD-GPR method varying
the number of snapshots. The number of POD modes is fixed to
$20$.\label{fig:sens_snaps}}
\end{figure}

\paragraph{Proper orthogonal decomposition with Gaussian process regression}
We exploit the collected database in order to build a kind of probabilistic
response surface to predict the resistance of new shapes. We remember the
starting set is composed by $100$ input-output pairs $\{(\mupar_i,
\out_i)\}_{i=1}^{100}$, where $\mupar$ is the geometrical parameters provided
to FFD and $\out$ is the resistance distribution over the deformed hull. Of
the entire set, we use the $80\%$ for train the POD-GPR framework and exploit
the remaining pairs to test our method.
Firstly we applied the POD on the snapshots matrix to reduce the dimension of
the output.
\begin{figure}[t]
\centering
\input{plots_arxiv/sing_POD.tikz}
\caption{Values of $\frac{\sigma_i}{\sigma_0}$ for the first 40 singular values, ordered from the highest to the lowest value.}
\label{sing_POD}
\end{figure}
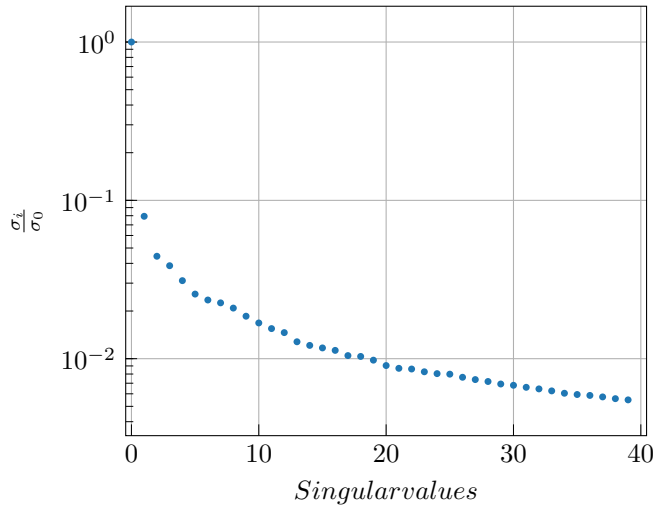
In this case the singular values extracted are reported in \autoref{sing_POD}
from an energetic perspective. The decay-rate is not very steep, probably due
to the discontinuous component for the VOF variable $\alpha$, which is directly
involved in the resistance computation.
\begin{figure}[t]
\centering
\input{plots_arxiv/err2modes.tikz}
\caption{Sensitivity analysis on the accuracy of POD-GPR method varying
the number of POD modes used. The number of snapshots is fixed to
$80$.\label{fig:sens_modes}}
\end{figure}
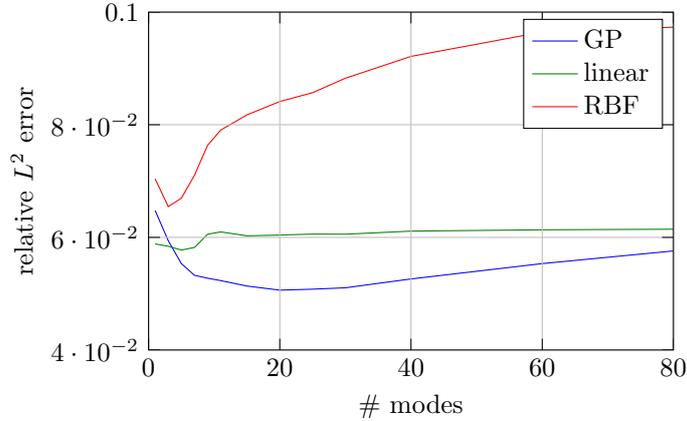
Despite this, the POD allows to remarkably reduce the dimension of the output,
simplifying the next phase.
We exploit the computed modal coefficient in order to optimize the GP, then query for
the new parametric solutions. 
To measure the accuracy, we propose in \autoref{fig:sens_modes} and
\autoref{fig:sens_snaps} two different sensitivity analysis varying the number
of POD modes used to span the reduced space and the number of snapshots to
train the method, respectively. For sake of completeness, we compare the results
with similar data-driven methodologies that involve, instead the GPR, other
interpolation techniques for modal coefficient approximation, as the linear
interpolation one or the radial basis function (RBF) one. 
We propose here the simplest RBF interpolation, but we make the reader aware
that better results can be achieved tuning the smoothness of RBF, producing a
non-interpolating RBF method.
For more details we refer~\cite{mythesis}.
The error refers to the mean relative error computed on the test dataset (of
dimension $20$), using the resistance distribution coming from DMD as truth
solution. The GPR method is able to reach the minimum error respect to the
other interpolation, resulting in a relative error near to $5\%$ adopting $20$
modes, but reducing its accuracy increasing the number of modes. This trend is
shown also by RBF error, that after an initial decreasing, becomes very large
for many modes.
Varying the number of snapshots (\autoref{fig:sens_snaps}), the difference
between RBF and GPR is even more evident. While the RBF reaches an error slightly
less than the $8\%$, the GPR is able to stay beyond the $6\%$ with $80$
snapshots. We note that we get the highest difference between the methods using
few snapshots: the GPR shows higher accuracy even with few samples and a
pretty constant trend for database with greater dimension.
\begin{figure}[t]
\centering
\includegraphics[width=.85\textwidth]{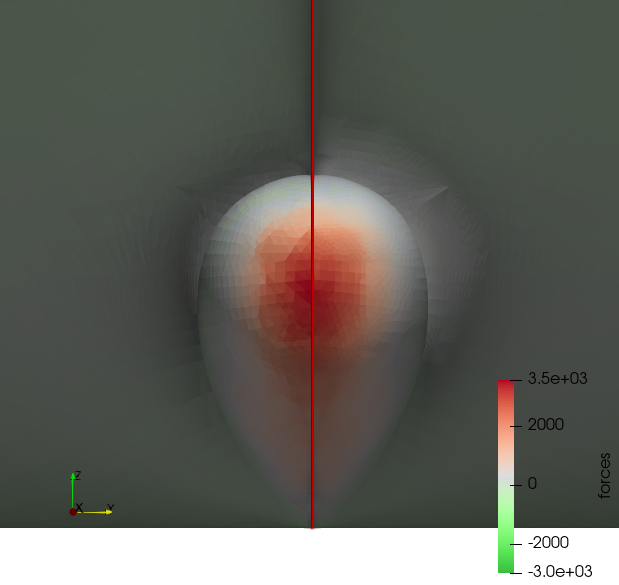}
\caption{Value of total resistance over the bulbous bow for the FOM (on the left) and for the ROM (on the right).}
\label{ROMFOMfrontal}
\end{figure}
Finally, we
conclude with a graphical visualization of the resistance distribution on (a
limited region\footnote{The bulbous bow is one of the region where the pressure
resistance is higher, and then difficult to predict.} of) the hull in
\autoref{ROMFOMfrontal}, comparing the ROM approximation with the FOM
validation. Even if the difference is notable, the reduced model can express
the main physics behaviour of the original model.
Regarding the computational cost reduction, the reduced model can approximate
the parametric solution only sampling an already defined distribution, and even
on a personal laptop it takes no more than few seconds, whereas the FV solver
takes $8$ hours, resulting in a very huge speedup.

\begin{figure}[t]
\centering
\input{plots_arxiv/opt1vr.tikz}
	\caption{Optimal shapes produced by $15$ different application of GA,
	in terms of resistance and volume as percentage respect to the original
	ship.}
\label{fig:opt}
\end{figure}
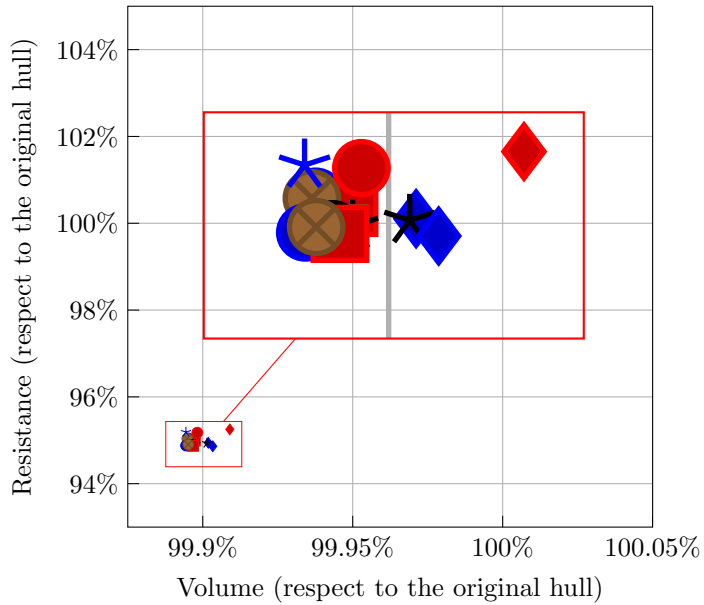

\paragraph{Genetic algorithm}
The goal of the entire pipeline is the minimization of the total resistance
(only in the direction of the flow). To ensure feasibility of the deformed
shape from the engineering viewpoint, we add a penalization on the hulls whose
volume is lower that $999\permil$ of the original hull. In other words, we
penalize the configurations that lead a volume decrease greater than
$\sfrac{1}{1000}$ with respect to the original volume. Our optimization problem
reads:
\begin{equation}
\underset{\mupar}{\text{min}}
\left\{
\begin{aligned}
	&\smallint_{\pdomain} \tau_x \rho - p n_x\quad & \text{if} \smallint_{\pdomain}\rho g h &\geq 0.999 \smallint_{\domain} \rho g h\\
	&\infty & \text{otherwise} &\\
\end{aligned}
\right.
	\label{eq:opt}
\end{equation}
where $\tau_x$ is the $x$ component of the (viscous and turbulent) tangential
stresses, $\rho$ is the density of the fluid (computed according to the VOF
model), $p$ is the pressure, $n_x$ is the $x$ component of the normal to the
surface, $g$ is the gravity acceleration and $h$ is the distance between the
surface and the waterline ($\smallint_\domain \rho g h$ results the volume of
the immersed hull using an hydrostatic approach).
To compute the objective function for a generic parameters, we need to perform
the FFD morphing then project the POD-GPR solutions over the deformed ship, in
order to numerically compute such integral. We clarify that with the reduced
order model returns the distribution of viscous and pressure forces over the
hull, that is $\tau \rho - p$ in $\pdomain$. As already mentioned, these methods
have a negligible computational cost, allowing us to optimize the shape in a
very efficient manner.
Despite its easiness of application, GA requires a good tuning of the
hyper-parameters to result successful. In this work, we applied the one point
crossover~\cite{eshelman1993real} for the mate procedure, while for the
mutation a Gaussian mutation~\cite{hinterding1995gaussian} with $\sigma = 0.1$
has been involved.  We set the mate and mutation
probability to $0.8$ and $0.2$, respectively. Moreover, we use an initial population
dimension $N_0 = 200$, reducing it to $N = 30$ during the evolution. The
stopping criteria in this case is the number of generations, which is set to
$15$. Robustness of this setting is proved in \autoref{fig:opt}, where $15$
different runs have been run and, for each run, the optimal shape is plotted
both in terms of resistance and volume. 
We can note in fact that all the runs have converged to the same fitness,
despite the stochastic component of the method itself, ensuring that the
hyper-parameters are set to fully explore the parameter space (and then globally
converge to the optimal point). The penalization we impose avoids the creation
of unfeasible deformations: the optimum of all the runs show a slight decreased
volume, but within the initial tolerance, while the resistance results decreased by
more than the $4\%$.

We specify that this is the optimum for the reduced model. In order to obtain an
accurate value, the optimal parameter can be plug in the pipeline and the
optimal shape is then validated using the full-order FV method. Additionally, this latter
can be insert in the snapshots database and used to enrich the precision of
the POD-GPR model. In our case, after the validation, the gain in term of
resistance is lower with respect to the ROM approximation, but reaching the $3.3\%$ it
results in a very good outcome in the engineering context.

%% file: plots_arxiv/eigs.tikz
\begin{tikzpicture}

\begin{axis}[
xlabel={Real part},
xmajorgrids,
xmin=-1.2, xmax=1.2,
xtick style={color=black},
ylabel={Imaginary part},
ymajorgrids,
ymin=-1.2, ymax=1.2,
ytick style={color=black},
 width=0.7\textwidth,
        height=0.7\textwidth,
]
\addplot [semithick, blue, mark=o, mark size=2, mark options={solid}, only marks]
table {%
0.999979368683168 0
-0.0220095980832188 0
0.841886605662969 0.541379001256738
0.841886605662969 -0.541379001256738
0.662514582624404 0
};
\addlegendentry{Eigenvalues}
\draw[red, dashed] (axis cs:0,0) circle [radius=82pt];

\end{axis}

\end{tikzpicture}

%% file: plots_arxiv/sing_POD.tikz
\begin{tikzpicture}
\definecolor{color0}{rgb}{0.12156862745098,0.466666666666667,0.705882352941177}

\begin{axis}[
tick align=inside,
tick pos=left,
x grid style={white!69.01960784313725!black},
xmin=-0.45, xmax=40.45,
ymode=log,
xtick style={color=black},
y grid style={white!69.01960784313725!black},
ytick style={color=black},
xlabel={$Singular values$},
ylabel={$\frac{\sigma_i}{\sigma_0}$},
ymajorgrids,
xmajorgrids,
]
\addplot [semithick, color0, mark=*, mark size=1, mark options={solid}, only marks]
table {%
0 1
1 0.079325930467808406
2 0.044358115629986053
3 0.038665325441975061
4 0.031111192974307089
5 0.025590394887916482
6 0.023476240196564562
7 0.022520008334246673
8 0.020875735706549343
9 0.018550386928347311
10 0.016815444149393573
11 0.015490448879560095
12 0.014612523163285023
13 0.012789273197667211
14 0.012132579941989621
15 0.01169516752639896
16 0.011284295856202132
17 0.010461538917619147
18 0.010331799257341947
19 0.0097886069612952344
20 0.0090498860312925155
21 0.0086927793248518948
22 0.0086060948843690449
23 0.0082665711404393393
24 0.0080467735057137592
25 0.0079857665099015308
26 0.0076324874796307096
27 0.0073803216577253026
28 0.0071780727900453587
29 0.0069202898565625072
30 0.0067931257282935371
31 0.0065959045404990159
32 0.0064450934388363967
33 0.0062607317144129371
34 0.0060495666660243391
35 0.0059394586011262832
36 0.0058615066247734446
37 0.0057356226055359375
38 0.0055919218952577444
39 0.0054936175506730533
};
\end{axis}

\end{tikzpicture}

%% file: plots_arxiv/err2modes.tikz
\begin{tikzpicture}

\definecolor{color1}{rgb}{0.75,0,0.75}
\definecolor{color2}{rgb}{0.75,0.75,0}
\definecolor{color0}{rgb}{0,0.75,0.75}

\begin{axis}[
axis on top,
legend cell align={left},
tick pos=both,
xmin=0, xmax=80,
xtick style={color=black},
ymin=0.04, ymax=0.1,
ytick style={color=black},
xlabel={\# modes},
ylabel={relative $L^2$ error},
ymajorgrids,
xmajorgrids,
width=.7\textwidth,
height=.5\textwidth,
]
\addplot [blue]
table {%
1 0.0647670728215381
3 0.0594114426723079
5 0.0553456394401349
7 0.0532667059959067
9 0.0527579392911333
11 0.0523239237504958
15 0.0513630872338408
20 0.0506309633501817
25 0.0508040216608267
30 0.0510512903515633
40 0.0526220437864548
60 0.0553470812622952
80 0.0575928650625103
};
\addlegendentry{GP}
\addplot [green!50.0!black]
table {%
1 0.0588433775747572
3 0.058424609774319
5 0.0577508923376143
7 0.0582170598286952
9 0.060558049317481
11 0.0609708436444476
15 0.0602791913752857
20 0.0604120970259428
25 0.060586459662219
30 0.0605634747339095
40 0.0611118263030476
60 0.0613459129629952
80 0.0614499113214857
};
\addlegendentry{linear}
\addplot [red]
table {%
1 0.0703925136041706
3 0.0654666824881262
5 0.0669587562137952
7 0.071035751636373
9 0.0763485993900294
11 0.0790590236081722
15 0.0817540323440976
20 0.084119749957877
25 0.0856928197649325
30 0.0882591140120198
40 0.092130711517646
60 0.096455359991419
80 0.0973185320922635
};
\addlegendentry{RBF}
\end{axis}

\end{tikzpicture}

%% file: plots_arxiv/err2snapshots.tikz
\begin{tikzpicture}

\definecolor{color1}{rgb}{0.75,0,0.75}
\definecolor{color0}{rgb}{0,0.75,0.75}

\begin{axis}[
axis on top,
legend cell align={left},
tick pos=both,
xmin=10, xmax=80,
xtick style={color=black},
ymin=0.04, ymax=0.1,
ytick style={color=black},
xlabel={\# snapshots},
ylabel={relative $L^2$ error},
ymajorgrids,
xmajorgrids,
height=.5\textwidth,
width=.7\textwidth,
]
\addplot [blue]
table {%
10 0.0641873803711008
15 0.0597589714254302
20 0.0582636993309825
30 0.0566366437175127
40 0.0553956167437841
60 0.0543279674642333
80 0.0535621058816762
};
\addlegendentry{GP}
\addplot [red]
table {%
10 0.0983259629012611
15 0.0865196690764611
20 0.0874302637727341
30 0.0837367872824286
40 0.0804280499747817
60 0.0779689702422913
80 0.0809664291099436
};
\addlegendentry{RBF}
\end{axis}

\end{tikzpicture}

%% file: plots_arxiv/opt1vr.tikz
\begin{tikzpicture}

\definecolor{color0}{rgb}{0,0.75,0.75}
\definecolor{color1}{rgb}{0.75,0,0.75}
\begin{scope}[spy using outlines={rectangle, magnification=5,
   width=5cm,height=3cm,connect spies}]

\begin{axis}[
tick align=inside,
tick pos=left,
x grid style={white!69.01960784313725!black},
xmajorgrids,
xmin=99.875, xmax=100.05,
xtick style={color=black},
y grid style={white!69.01960784313725!black},
ymajorgrids,
xticklabels={99.9\%, 99.95\%, 100\%, 100.05\%},
xtick={99.9, 99.95, 100, 100.05},
yticklabels={94\%, 96\%, 98\%, 100\%, 102\%, 104\%},
ytick={94, 96, 98, 100, 102, 104},
ylabel={Resistance (respect to the original hull)},
xlabel={Volume (respect to the original hull)},
ymin=93, ymax=105,
ytick style={color=black},
height=.7\textwidth,
width=.7\textwidth,
]
\addplot table {
	99.8951261645666 95.0442026420136
};
\addplot table {
	99.8973624826236 94.9918984515746
};
\addplot table {
	99.8948803574267 95.0369530136232
};
\addplot table {
	99.8975976016157 94.9073793940933
};
\addplot table {
	99.9018290926176 94.942429275442
};
\addplot table {
	99.8981791959246 95.1756292931723
};
\addplot table {
	99.8967438188275 94.8760912945524
};
\addplot table {
	99.896269061514 94.8942875059163
};
\addplot table {
	99.8943902435377 95.1895275591046
};
\addplot table {
	99.9090231639592 95.2515431274815
};
\addplot table {
	99.894439403734 94.8784143552404
};
\addplot table {
	99.8967490627382 94.8726429192259
};
\addplot table {
	99.8951648301315 94.9038124033823
};
\addplot table {
	99.9014021598079 94.9346594639952
};
\addplot table {
	99.9033326953859 94.8625142993142
};
\end{axis}
\spy [red] on (1,1.1) in node (zoom) [left] at ([xshift=6cm,yshift=4cm]0, 0);
\end{scope}

\end{tikzpicture}

%% file: sections/conclusion.tex
\section{Conclusion and future perspectives}
\label{sec:conclusion}

In this work, we propose a complete computational framework for shape optimization problems.
To overcome the computational barrier, the dynamic mode decomposition (DMD) and
the proper orthogonal decomposition with Gaussian process regression (POD-GPR)
are involved. This pipeline aims at reducing the number of high-fidelity
simulation needed to converge to the optimal shape, making its application very
useful in all the context where the performance evaluation of the studied
object results computationally expensive.
We applied such framework to an industrial shape optimization problem,
minimizing the total resistance of a cruise ship advancing in calm water.
Exploitation of ROM techniques drastically reduces the overall time, and even if
the accuracy of the reduced model is decreased (with respect to the full-order one) the
final outcome presents a remarkable reduction of the resistance ($3.3\%$).

Future developments regarding this integrated methodology may interest the
extension to constraint optimization problems, the involvement of machine
learning techniques in the optimization procedure, or a greedy approach that
enriches the reduced order model by adding iteratively the approximated optimal
shape.